\documentclass[12pt]{amsart}
\usepackage{fancyhdr} 
\usepackage{setspace}
\usepackage{titlesec}
\usepackage{enumitem} 
\usepackage{amsmath}
\usepackage{amsthm}
\usepackage{amssymb}
\usepackage{multicol}
\usepackage{enumerate,pifont,mathrsfs,geometry}
\usepackage{tikz}
\usepackage{tikz-cd}
\usepackage{mathtools}
\usepackage[numbers]{natbib}
\usepackage{float}
\usepackage{pgfplots}
\usepackage{geometry}
\usepackage{xcolor}
\usepackage{hyperref}
\hypersetup{
    colorlinks=true,
    linkcolor=red,
    citecolor=blue,
    urlcolor=blue
}

\pgfplotsset{compat=newest}
\usepgfplotslibrary{fillbetween}
\allowdisplaybreaks

\newtheorem{theorem}{Theorem}[section]
\newtheorem{corollary}[theorem]{Corollary}

\theoremstyle{definition}
\newtheorem{definition}[theorem]{Definition}

\def\subsection{\def\@secnumfont{\bfseries}\@startsection{subsection}{1}%
  {\parindent}{.5\linespacing\@plus.7\linespacing}{-.5em}%
  {\normalfont\bfseries}}
\def\pts#1{\hfill({\em #1 point\ifnum #1>1 s\fi})\ifnum#1>9\else\phantom{\em 1{\ifnum#1>1\else s\fi}}\fi}

\makeatother
 \titleformat{\section}
  {\normalfont\scshape \filcenter}{\thesection}{.5em}{}

\newcommand{\R}{\mathbb{R}}

\newcommand{\ob}{\mathfrak{o}}

\usepackage{amsmath,amssymb,amsthm, mathrsfs, mathtools}
\usepackage{amscd,mathtools}
\usepackage{fullpage}
\usepackage{graphicx}
\usepackage{color,xcolor}
\usepackage{enumerate}
\usepackage{xspace}
\usepackage{tipa}
\usepackage{stmaryrd}
  
\usepackage[utf8]{inputenc}
\usepackage[T1]{fontenc} \usepackage{imakeidx}
\makeindex[columns=3, title=Alphabetical Index, intoc]

\usepackage{tikz}
\usetikzlibrary{shapes.geometric}
\usepackage{caption, subcaption}
\usepackage{tikz-cd}
\usetikzlibrary{calc,decorations.pathreplacing}
\usepackage{tikz}
\usetikzlibrary{calc}
\usetikzlibrary{arrows}
\usetikzlibrary{decorations.pathreplacing}
\usetikzlibrary{intersections}
\usetikzlibrary{decorations.pathmorphing}

\usepackage{mathrsfs}
\usetikzlibrary{matrix}
\usetikzlibrary{positioning}
\DeclareMathAlphabet{\pazocal}{OMS}{zplm}{m}{n}
\tikzset{>=stealth}

\usepackage{hyperref}

  \theoremstyle{definition}

  \newtheorem*{claim*}{Claim}

  \newtheorem*{question*}{Question}
  \newtheorem*{answer*}{Answer}
  \newtheorem*{application*}{Application}

  \theoremstyle{remark}
  \newtheorem*{remark*}{Remark}
  
  \newcommand{\param}{{\mathchoice{\mkern1mu\mbox{\raise2.2pt\hbox{$
  \centerdot$}}
  \mkern1mu}{\mkern1mu\mbox{\raise2.2pt\hbox{$\centerdot$}}\mkern1mu}{
  \mkern1.5mu\centerdot\mkern1.5mu}{\mkern1.5mu\centerdot\mkern1.5mu}}}

\DeclarePairedDelimiterX{\norm}[1]{\lvert}{\rvert}{#1}
\DeclarePairedDelimiterX{\Norm}[1]{\lVert}{\rVert}{#1}

\title{The Curtain Model is Not a Quasi-Isometry Invariant of CAT(0) Spaces}
\author{Elliott Vest}
\begin{document}

\maketitle
\begin{abstract}
    Petyt-Spriano-Zalloum recently developed the notion of a \textit{curtain model}, which is a hyperbolic space associated to any CAT(0) space.  It plays a similar role for CAT(0) spaces that curve graphs do for mapping class groups of finite-type surfaces.  Those authors asked whether this curtain model is a quasi-isometry invariant, namely if quasi-isometric CAT(0) spaces have quasi-isometric curtain models.  In this short note, we provide an explicit example answering this question in the negative.
\end{abstract}

\section{Introduction}

 In \cite{PSZ22}, Petyt-Spriano-Zalloum introduced a combinatorial tool called a \textit{curtain} that serves as an analogue in the CAT(0) setting of a hyperplane from the theory of CAT(0) cube complexes. Building off of ``hyperplane-separation" metrics introduced by Genevois \cite{Gen20b}, the authors utilize curtains in a CAT(0) space $X$ to build the \textit{curtain model} --- a hyperbolic space which effectively collapses the ``flat" parts of the CAT(0) space. This ``coning off" is by design, and gives rise to many similarities between a CAT(0) space and its curtain model that parallel the relationship between mapping class groups and their curve graphs (See \cite{Zal23}).

Petyt-Spriano-Zalloum asked in \cite{PSZ22} if a quasi-isometry between CAT(0) spaces always induces a quasi-isometry between their corresponding curtain models. We answer this question in the negative. For a CAT(0) space $X$, we denote $\widehat X$ to be its curtain model.

\begin{theorem}\label{nonQI}
    There exists a CAT(0) space $X$ and a self quasi-isometry $\phi: X\longrightarrow X$ such that $\phi$ does not descend to a quasi-isometry for $\widehat X$. Further, there exists two quasi-isometric CAT(0) spaces $W, Z$ whose curtain models $\widehat W, \widehat Z$ are not quasi-isometric.
\end{theorem}

Our example is based on an example due to Cashen \cite{Cashen+2016}, which he used to show that quasi-isometries of CAT(0) spaces need not induce homeomorphisms of their contracting boundaries when equipped with the Gromov product topology. Thus, it also follows that we get an analogous result for the curtain models of CAT(0) spaces.

\begin{corollary}
    There exist quasi-isometric CAT(0) spaces $W,Z$ whose curtain models have non-homeomorphic Gromov boundaries.
\end{corollary}

\noindent \textbf{Acknowledgments:} Harry Petyt independently discovered this example.  I would like to thank him for his useful comments on an earlier draft of this article and for encouraging me to write it up. Also, the warmest of thanks goes to Matthew Gentry Durham for his constructive feedback on an earlier draft of this paper.
\section{Background}

We now give a small summary of definitions imported from \cite{PSZ22}. For background of CAT(0) spaces, we refer the reader to \cite[II.1]{BH99}. The following is the required background to define the \textit{curtain model} (Definition \ref{Curtain Model}). We always assume $X$ is a CAT(0) space.

 \begin{definition}[Curtain, Pole]\label{curtain} Let $X$ be a CAT(0) space and let $\alpha:I\rightarrow X $ be a geodesic. For any number $r$ such that $[r-\frac{1}{2}, r+\frac{1}{2}]$ in in the interior of $I$, the \textit{curtain dual to $\alpha$} at $r$ is 
 \begin{center}
 $h=h_\alpha=h_{\alpha,r}=\pi^{-1}_\alpha(\alpha[r-\frac{1}{2}, r+\frac{1}{2}])$
 \end{center}
where $\pi_\alpha$ is the closest point projection to $\alpha$. We call the segment $\alpha[r-\frac{1}{2}, r+\frac{1}{2}]$ the \textit{pole} of the curtain which we denote as $P$ when needed.
\end{definition}

\begin{definition}[Chain, Separates]
    A curtain $h$ \textit{separates} sets $A,B \subset X$ if $A \subset h^-$ and $B \subset h^+$. A set $\{h_i\}$ is a \textit{chain} if each of the $h_i$ are disjoint and $h_i$ separates $h_{i-1}$ and $h_{i+1}$ for all $i$. We say a chain $\{h_i\}$ \textit{separates} sets $A,B \subset X$ if each $h_i$ separates $A$ and $B$.
\end{definition}

\begin{definition}[$L$-separated, $L$-chain] \label{L chain}  Let $L \in \mathbb{N}$. Disjoint curtains $h$ and $h^{\prime}$ are said to be \textit{$L$-separated} if every chain meeting both $h$ and $h^{\prime}$ has cardinality at most $L$. Two disjoint curtains are said to be \textit{separated} if they are $L$-separated for some $L$. If $c$ is a chain of curtains such that each pair is $L$-separated, then we refer to $c$ as an \textit{$L$-chain}. 
\end{definition}

\begin{definition}[$L$-metric] \label{Lmetric} Denote $X_L$ for the metric space $(X, d_L)$, where $d_L$ is the metric defined as \begin{center}
     $\mathrm{d}_L(x, y)=1+\max \{|c|: c$ is an $L$-chain separating $x$ from $y\}$
\end{center}
with $d_L(x,x) = 0$. Note that, by Remark 2.16 in \cite{PSZ22}, we have that for any $x,y \in X$, it follows that $d_L(x,y) < 1 +d(x,y).$
\end{definition}

\begin{definition}[Curtain Model] \label{Curtain Model} Fix a sequence of number $\lambda_L \in (0,1)$ such that $$\sum_{L=1}^\infty \lambda_L < \sum_{L=1}^\infty L\lambda_L < \sum_{L=1}^\infty L^2\lambda_L < \infty$$
We consider the space $(X, \hat{d})$, where the distance between two points $x,y\in X$ is defined by \begin{center}
    $\hat{d}(x,y) = \displaystyle \sum_{L=1}^\infty \lambda_L d_L(x,y)$
\end{center}
and $d_L$ is the $L$-metric defined in Definition \ref{Lmetric}. We call $(X, \hat{d})$ the \textit{curtain model} of $X$ and denote it as $\widehat{X}$.
\end{definition}

Both of the following definitions will also help in the construction of the counterexample.

\begin{definition}[Angles in CAT(0) spaces, Section II.3.1 in \cite{BH99}]
Let $X$ be a CAT(0) space and let $\alpha:[0, a] \rightarrow X$ and $\alpha^{\prime}:\left[0, a^{\prime}\right] \rightarrow X$ be two geodesic paths issuing from the same point $\alpha(0)=\alpha^{\prime}(0)$. Then the comparison angle $\angle_{\mathbb{E}}\left(\alpha(t), \alpha^{\prime}\left(t^{\prime}\right)\right)$ is a non-decreasing function of both $t, t^{\prime} \geq 0$, and the \textit{Alexandrov angle} $\angle\left(\alpha, \alpha^{\prime}\right)$ is equal to
$$
\lim _{t, t^{\prime} \rightarrow 0} \angle_{\mathbb{E}}\left(\alpha(t), \alpha^{\prime}\left(t^{\prime}\right)\right)=\lim _{t \rightarrow 0} \angle_{\mathbb{E}}\left(\alpha(t), \alpha^{\prime}(t)\right) .
$$
Hence, we define:
$$
\angle\left(\alpha, \alpha^{\prime}\right)=\lim _{t \rightarrow 0} 2 \arcsin \frac{1}{2 t} d\left(\alpha(t), \alpha^{\prime}(t)\right) .
$$
\end{definition}

\begin{definition}[Strongly Contracting]
A geodesic $\alpha$ is $D$-strongly contracting if for any ball B disjoint from $\alpha$ we have $diam(\pi_\alpha(B)) \leq D$, where $\pi_\alpha$ is the closest point projection to $\alpha$.
    
\end{definition}

\section{The Counterexample}

The following counterexample was used in \cite{Cashen+2016} to show that two quasi-isometric CAT(0) spaces can have contracting boundaries of different homeomorphism type when equipped with the Gromov product topology. We first introduce this space and its curtain model.\vspace{.25cm}

\noindent \textbf{3.1 The Infinite Parking Lot and its Curtain Model.} \vspace{.25cm}

\begin{figure}[h]
    \centering
    \begin{tikzpicture} [scale=.4]

   \draw [->, gray] (0,1.65) to (0, 5.5);
   \draw [->, gray] (-1.8,0) to (-5.65,0);
   \draw [->, gray] (0,-2) to (0,-5.85);
   \draw [->, gray] (2.15,0) to (6,0);
   \draw [ gray] (1.45,0.1) to (2.15,0.1);


\fill[green, opacity=0.2] (0,-2) .. controls(.72,-2) and (1.85,-1.7) .. (2.15, 0) -- (5.7,0) arc (0:-90:5.7)  -- cycle;

\fill[yellow, opacity=0.3] (0,-2) .. controls(-.42,-2) and (-1.75,-1.7) .. (-1.8, 0) -- (-5.6,0) arc (180:270:5.6) -- cycle;

\fill[orange, opacity=.4] (-1.8,0) .. controls(-1.8,.95) and (-0.8,1.65) .. (0,1.6) -- (0, 5.55) arc (90:180:5.55)  -- cycle;

\fill[red, opacity=0.5] (0,1.6) .. controls(.1,1.65) and (.3,1.6) .. (.41,1.5) -- (0.41, 5.55) arc (85.9:90:5.55)  -- cycle;

\fill[red, opacity=0.5] (.41,1.5) .. controls(1.4,1.2) and (1.4,.2) .. (1.45,0.1) -- (2.15,0.1) .. controls(2.2,1.6) and (.8,2.2) ..(.41, 2.3)  -- cycle;

\draw [ gray, dotted] (2.15,0.1) to (6,0.1);


\fill[blue, opacity=0.1] (.41, 2.3)  .. controls(.8,2.2) and (2.2,1.6) ..(2.15,0)  -- (5.7,0) arc (0:84.1:5.9) -- cycle;



\draw (-9+.5,7+1) arc (0:90:1);
\draw (-10+.5,8+1) to (-10+.5,10+1);
\draw (-9+.5,7+1) to (-7+.5,7+1);

\fill[red, opacity=0.5] (-9+.5,7+1) arc (0:90:1)  -- (-9.5,11) arc (90:0:3) -- cycle;


\draw (-9+4+.5,7+1) arc (0:90:1);
\draw (-10+4+.5,8+1) to (-10+4+.5,10+1);
\draw (-9+4+.5,7+1) to (-7+4+.5,7+1);

\fill[orange, opacity=0.4] (-9+4+.5,7+1) arc (0:90:1) -- (-5.5,11) arc (90:0:3) -- cycle;


\draw (-9+8+.5,7+1) arc (0:90:1);
\draw (-10+8+.5,8+1) to (-10+8+.5,10+1);
\draw (-9+8+.5,7+1) to (-7+8+.5,7+1);

\fill[yellow, opacity=0.3] (-9+8+.5,7+1) arc (0:90:1) -- (-10+8+.5,11) arc (90:0:3) -- cycle;


\draw (-9+12+.5,7+1) arc (0:90:1);
\draw (-10+12+.5,8+1) to (-10+12+.5,10+1);
\draw (-9+12+.5,7+1) to (-7+12+.5,7+1);

\fill[green, opacity=0.2] (-9+12+.5,7+1) arc (0:90:1) -- (-10+12+.5,11) arc (90:0:3) -- cycle;


\draw (-9+16+.5,7+1) arc (0:90:1);
\draw (-10+16+.5,8+1) to (-10+16+.5,10+1);
\draw (-9+16+.5,7+1) to (-7+16+.5,7+1);

\fill[blue, opacity=0.1] (-9+16+.5,7+1) arc (0:90:1) -- (-10+16+.5,11) arc (90:0:3) -- cycle;


\node [below,] at (-7.7,7.8) {\Large$X_{i-1}$};
\node [below,] at (-7.7+4,7.8) {\Large$X_{i}$};
\node [below,] at (-7.7+8,7.8) {\Large$X_{i+1}$};
\node [below,] at (-7.7+12,7.8) {\Large$X_{i+2}$};
\node [below,] at (-7.7+16,7.8) {\Large$X_{i+3}$};

 \draw[thick, <->]
    plot[domain=720:1160, samples=144] ({\x}: {.2 * \x/100});

\end{tikzpicture}
    \caption{One level of the infinite parking lot $X = \cup_i X_i\slash \sim$. The space would continue to spiral upward and downward. Notice, since we are viewing the space from a birds eye view, $X_{i-1}$ gets shadowed by $X_{i+3}$. This is due to $X_{i+3}$ forming the ``next level" of the infinite parking lot.}
    \label{fig:parking lot}
\end{figure}
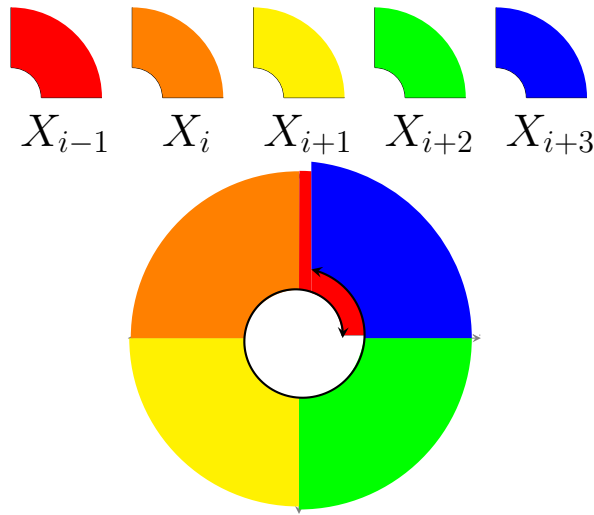

Let $Y$ be $\R^2$ with a disc of radius one centered at the origin removed. Denote $X$ as the universal cover of $Y$. We can view $X$ in the following way: Take $X_i$ to be a quarter flat with the quarter disc centered at the origin removed. Then $X = \cup_i X_i\slash \sim$ where $\sim$ denotes gluing the $y$-axis of $X_i$ to the $x$-axis of $X_{i+1}$ for all $i\in \mathbb{Z}$. One informally calls $X$ the ``infinite parking lot" as it can be viewed as a collection of quarter flats glued together that are spiraling up and down, giving the ``infinite levels" of a parking lot. See Figure \ref{fig:parking lot}.

  $X$ is indeed a CAT(0) space since it is a gluing of CAT(0) spaces along single geodesic lines. The result of this space is that a half flat with a half disc of radius one removed at the origin can be isometrically embedded into each $X_i\cup X_{i+1}\slash \sim$. In fact, we can spiral up any $\theta$ amount and get the same isometry of the half flat with a half disc removed at the origin. Parameterize $X$ via its natural polar coordinates $\R \times [1,\infty) $, and define \textit{the spiral} to be the line $\mathbb R \times \{1\}$. We now explain why $X$'s curtain model $\widehat X$ is a quasi-line.

\begin{figure}[h]
    \centering
    \begin{tikzpicture} [rotate=-90, scale=.5]

\fill[yellow, opacity=0.2] (0,-2) .. controls(-.42,-2) and (-1.75,-1.7) .. (-1.8, 0) -- (-7.6,0) arc (180:270:7.6) -- cycle;

\fill[orange, opacity=.3] (-1.8,0) .. controls(-1.8,.95) and (-0.8,1.65) .. (0,1.5) -- (0, 7.55) arc (90:180:7.55)  -- cycle;

\fill[white] (0.1,2.67) -- (-0.43,2.67) -- (-.43,7.6) -- (0.1,7.6) -- cycle;

 \node [below,] at (-3,4) {\Large$X_{i}$};
 \node [below,] at (-3,-4) {\Large$X_{i+1}$};

 \fill[blue, opacity=0.3] (-6.25,-4.3) -- (-5.75,-5) -- (-5.75,4.9) -- (-6.25,4.25) -- cycle;

  \fill[blue, opacity=0.3] (-6.25+1.25,-5.7) -- (-5.75+1.25,-6.1) -- (-5.75+1.25,6.05) -- (-6.25+1.25,5.65) -- cycle;

   \node [right, blue, opacity = 0.7] at (-5.75+1.0,6.05) {$h_1$};

   \fill[blue, opacity=0.3] (-6.25+2.5,-6.6) -- (-5.75+2.5,-6.85) -- (-5.75+2.5,6.8) -- (-6.25+2.5,6.55) -- cycle;

   \node [right, blue, opacity = 0.7] at (-5.75+2.3,6.8)  {$h_2$};

   \draw [ gray] (0,1.5) to (0, 2.67);
   \draw[dotted,gray] (0,2.67) to (0,7.85);
   \draw [->, thick, gray] (-1.8,0) to (-7.65,0);
    \node [above,] at (-7.65,0) {\Large$\gamma$};
   \draw [ gray] (0,-2.1) to (0,-7.85);
   \draw[thick, ->]
    plot[domain=0:820, samples=144] ({\x}: {1.2 * \x/360});
\end{tikzpicture}
    \caption{Since $\gamma$ can always be seen as in the middle of a half flat, its dual curtains will behave like curtains in a half flat. This implies that any pair of curtains dual to $\gamma$ will not be $L$-separated for any $L$.}
    \label{fig:gamma}
\end{figure}
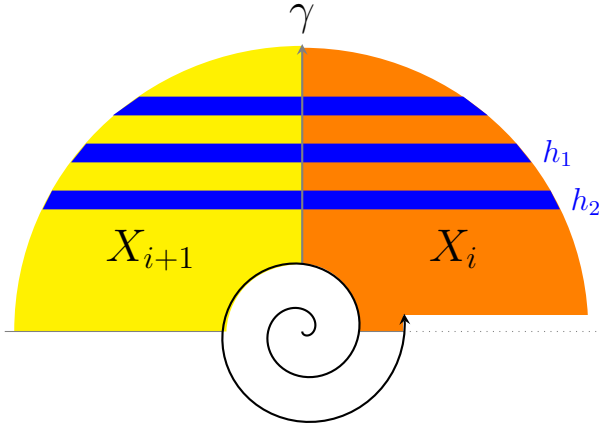

Take any geodesic ray $\gamma$ such that $\gamma(0)$ is on the spiral and the Alexandrov angle between $\gamma$ and the spiral is $\frac{\pi}{2}$. Up to an isometric rotation of $X$ by some $\theta$ along the spiral, $\gamma$ is the $y$-axis of some $X_i$. Since $\gamma$ is the $y$-axis of some isometrically embedded half flat (with a half disc removed), all curtains dual to $\gamma$ will stay in its half flat, $X_i\cup X_{i+1}\slash \sim$. As seen in Figure \ref{fig:gamma}, if $h_1, h_2$ are two disjoint curtains dual to $\gamma$, then $h_1,h_2$ will be two parallel, infinitely long strips of width one in $X_i\cup X_{i+1}\slash \sim$. All curtains dual to the $x$-axis of $X_i$ will meet $h_1$ and $h_2$, which means $h_1$ and $h_2$ are not $L$-separated for any $L$. The same is true for any two disjoint curtains dual to $\gamma$. Also, by Lemma 2.21 in \cite{PSZ22}, the max $L$-chain that can cross $\gamma$ is bounded above by $4L+10$. Thus, the diameter of $\gamma$ is $$\widehat{diam}(\gamma) = \sum_{L=1}^\infty \lambda_Ldiam_L(\gamma) \leq \sum_{L=1}^\infty \lambda_L(4L+10) < \infty.$$ This is true for any geodesic ray that starts at the spiral and whose Alexandrov angle with the spiral is $\frac{\pi}{2}$. In particular, if we denote the spiral as $\alpha$, then for any $x \in X$, $\hat d (x, \pi_\alpha(x)) \leq 4L+10$.

Now, fix some origin $\ob \in X$ on $\alpha$, and let $\alpha^+$ denote the positive spiral direction and $\alpha^-$ the negative spiral direction emanating from $\ob$. Both directions are $\pi$-strongly contracting as balls disjoint from the axis can only project to half of the circumference of one of the circles in the spiral. By \cite[Theorem 4.2]{PSZ22}, there exists an infinite $L$-chain dual to $\alpha^+$ for some $L$ (similarly for $\alpha^-$). Thus, in the curtain model $\widehat X$, the diameters of $\alpha^+$ and $\alpha^-$ will both be unbounded. By \cite[Proposition 9.5]{PSZ22}, both $\alpha^+$ and $\alpha^-$ are unparameterized quasi-geodesics in $\widehat X$. This concludes $\alpha$ is a quasi-line in $\widehat X$. Since for any $x \in X$, $\hat d(x, \pi_{\alpha}(x))\leq 4L+10$, this yields that $\widehat X$ is a quasi-line. \vspace{.25cm}

\noindent \textbf{3.2 A Self Quasi-Isometry Does Not Induce a Quasi-Isometry of Curtain Models.} \vspace{.25cm}

For some $\ob \in X$ on the spiral, denote the points of $X$ by $(\theta, r)$, where $\theta$ is the angle traveled around the spiral starting at $\ob$, and $r$ is the ``radius" distance away from the spiral. Consider the points $(i, 2^i)$ and $(0, 2^i)$  for all $i \in \mathbb{N}.$ Through a variation of the logarithmic spiral quasi-isometry of the Euclidean plane \begin{align*}
    \phi:X &\longrightarrow X\\
      (t,r)&\longmapsto (t-\log_2(r),r),
\end{align*} we see that $\phi\left( (i, 2^i)\right) = (0, 2^i).$ However, in the curtain model $\widehat{X}$, $\{(0,2^i)\}_i$ represents a quasi-point, and $\{(i,2^i)\}_i$ represents a quasi-line. This means that the self-quasi-isometry $\phi$ will not descend to a quasi-isometry for $\widehat{X}$.\vspace{.25cm}

\noindent \textbf{3.3 Upgrading to a Counterexample for Quasi-Isometric Invariance.}\vspace{.25cm}

\begin{figure}[ht]
    \centering
    \begin{tikzpicture}[scale=.5]
    \draw[thick, <->] (-7,0) to (7,0);
    \node[left]  at (-7,0.2) {\huge$\widehat{Z}$};

    \draw[thick, <->] (-7,3) to (7,3);
    \node[left]  at (-7,3.2) {\huge$\widehat{W}$};

    \node[red]  at (0,3) {\Large$\times$};
    \node[red]  at (1,3) {\Large$\times$};
    \node[red]  at (2,3) {\Large$\times$};
    \node[red]  at (3,3) {\Large$\times$};
    \node[red]  at (6,3) {\Large$\times$};
    \node[ below]  at (0,2.9) {\Large$0$};
    \node[ below]  at (1,2.9) {\Large$1$};
    \node[ below]  at (2,2.9) {\Large$2$};
    \node[ below]  at (3,2.9) {\Large$3$};
    \node[ below]  at (6,2.9) {\Large$i$};
    \node[ below]  at (4.5,2.9) {\Large$\cdots$};

    \node[ below]  at (0,-.1) {\Large$0$};
    \node[ below]  at (1,-.1) {\Large$1$};
    \node[ below]  at (2,-.1) {\Large$2$};
    \node[ below]  at (3,-.1) {\Large$3$};
    \node[ below]  at (6,-.1) {\Large$i$};
    \node  at (1,0) {\tiny$|$};
    \node  at (2,0) {\tiny$|$};
    \node  at (3,0) {\tiny$|$};
    \node  at (6,0) {\tiny$|$};
    \node[ below]  at (4.5,-.1) {\Large$\cdots$};

    \node[red]  at (0,0) {\Large$\times$};
    \node[red]  at (-.2,.2) {\Large$\times$};
    \node[red]  at (.1,.3) {\Large$\times$};
    \node[red]  at (-.3,-.2) {\Large$\times$};
    \node[red]  at (.4,.1) {\Large$\times$};
    \node[red]  at (-.1,.5) {\Large$\times$};
\end{tikzpicture}

    \caption{In this picture, the red marks represent the starting points of the geodesic rays that were glued to $X$. Notice that for $\widehat W$, there is a geodesic ray starting at each natural number. On the other hand, the geodesic rays in $\widehat Z$ are all clumped at the origin. The curtain models of $W$ and $Z$ both crunch the flatness of the parking lot, but the geodesic rays still have recognizable distance at their starting points in $\widehat W$. Thus, the two curtain models are not quasi-isometric.}
    \label{fig:Z and W}
\end{figure}
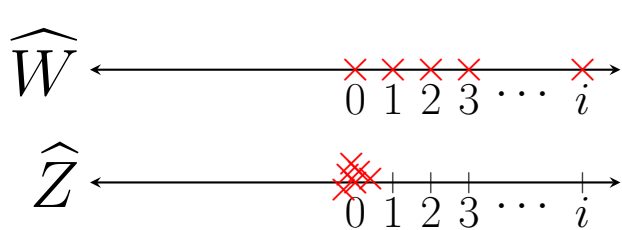

Now, following the same vein as \cite{Cashen+2016}, we construct two quasi-isometric CAT(0) spaces whose curtain models are not quasi-isometric. Construct the space $W$ by gluing a geodesic ray $\gamma_i$ to $X$ at each $(i, 2^i)$ point. Similarly, construct the space $Z$ by gluing a geodesic ray $\gamma_i'$ to $X$ at each $(0, 2^i)$ point. These spaces are quasi-isometric via the quasi-isometry \begin{align*}
    \overline \phi:W &\longrightarrow Z\\
      (t,r)&\longmapsto (t-\log_2(r),r)\\
      \gamma_i&\longmapsto \gamma_i'.
\end{align*} However, the curtain models will not be quasi-isometric. See Figure \ref{fig:Z and W}. Indeed, as $\{(0,2^i)\}_i$ is a quasi-point in $\widehat Z$,  each of the geodesic rays in $\{\gamma_i'\}_i$ emanate from a point which is within bounded distance of $\ob$ on the quasi-line $\widehat X$. Thus, $\widehat Z$ is quasi-isometric to an infinite wedge of rays. On the other hand, $\{(i,2^i)\}_i$ represents some sub-quasi-line in $\widehat X$, so the geodesic rays $\{\gamma_i\}_i$ have starting points at increasing distance away from $\ob$ in $\widehat X$ as $i$ increases. So, $\widehat W$ is quasi-isometric to $\mathbb R$ with a ray attached to each positive integer. These two spaces are not quasi-isometric.

The same logic can also apply to show $\widehat W$ and $\widehat Z$ have Gromov boundaries of different homeomorphism type. The sequence $\{\gamma_i\}_i$ in the Gromov boundary of $\widehat W$ converges to $\alpha^+$. No such converging sequence exists in $\widehat Z$. This proves the two Gromov boundaries for $\widehat W, \widehat Z$ are not homeomorphic.

\bibliography{QI-Invariance}{}
\bibliographystyle{alpha}

\end{document}